\documentclass[12pt]{amsart}

\usepackage{amssymb, amscd, txfonts}
\usepackage{graphicx}

\numberwithin{equation}{section}

\newtheorem{thm}{Theorem}[section]
\newtheorem{prop}[thm]{Proposition}
\newtheorem{lem}[thm]{Lemma}
\newtheorem{cor}[thm]{Corollary}

\theoremstyle{definition}
\newtheorem{defn}[thm]{Definition}

\theoremstyle{remark}
\newtheorem{rem}[thm]{Remark}
\newtheorem{exmp}[thm]{Example}

\renewcommand{\hom}{\operatorname{Hom}}
\renewcommand{\ker}{\operatorname{Ker}}

\newcommand{\Z}{\mathbb{Z}}

\newcommand{\R}{\mathbb{R}}
\newcommand{\C}{\mathbb{C}}
\newcommand{\F}{\mathbb{F}}

\DeclareMathOperator{\Ad}{Ad}
\DeclareMathOperator{\aut}{Aut}

\DeclareMathOperator{\im}{Im}
\DeclareMathOperator{\out}{Out}
\DeclareMathOperator{\sgn}{sgn}
\DeclareMathOperator{\tr}{tr}
\DeclareMathOperator{\tor}{Tor}
\DeclareMathOperator{\vol}{Vol}

\DeclareMathOperator{\eul}{Eul}

\begin{document}

\title[Symmetry of Reidemeister torsion on $SU_2$-representation spaces of knots]{Symmetry of Reidemeister torsion on $SU_2$-representation spaces of knots}
\author[T.~Kitayama]{Takahiro KITAYAMA}
\address{Graduate~School~of~Mathematical~Sciences, the~University~of~Tokyo, 3-8-1~Komaba, Meguro-ku, Tokyo 153-8914, Japan}
\email{kitayama@ms.u-tokyo.ac.jp}

\begin{abstract}
We study two sorts of actions on the space of conjugacy classes of irreducible $SU_2$-representations of a knot group.
One of them is an involution which comes from the algebraic structure of $SU_2$ and the other is the action by the outer automorphism group of the knot group.
In particular, we consider them on a $1$-dimensional smooth part of the space, which is canonically oriented and metrized via a Reidemeister torsion volume form.
As an application we show that the Reidemeister torsion function on the $1$-dimensional subspace has symmetry about the metrization.
\end{abstract}

\maketitle

\section{Introduction}

This work was intended as an attempt to globally describe Reidemeister torsion as a function on the $SU_2$-representation space of a knot group and show symmetry of the function by considering natural actions on the representation space.
Reidemeister torsion of a knot associated to a linear representation of the knot group is known to coincide with the twisted Alexander invariant associated to the representation.
Therefore this study also deduces symmetry of the twisted Alexander function.
See \cite{KL}, \cite{Ki}, \cite{L} and \cite{Wa} for the definition of twisted Alexander invariants and the relation with Reidemeister torsion.

Let $K$ be an oriented knot in an oriented rational homology $3$-sphere and $E$ the complement of an open tubular neighborhood of $K$.
We consider a $1$-dimensional smooth part $\mathcal{R}$ of the space of conjugacy classes of irreducible $SU_2$-representations of $\pi_1 E$.
(defined in Section \ref{subsec_2_1}.)
For $K$ in $S^3$, $\mathcal{R}$ is the set of representations which are called regular.
In this case Heusener \cite{H} showed that $\mathcal{R}$ carries a canonical orientation and Dubois \cite{D1}, \cite{D2} established a canonical volume form $\tau$ on $\mathcal{R}$ which induces Heusener's orientation via non-acyclic Reidemeister torsion. 
The idea of regarding non-acyclic Reidemeister torsion as volume forms on representation spaces was first considered by Witten.
In \cite{Wi} he obtained a remarkable formula to compute the volume of the moduli spaces of representations of fundamental groups of surfaces in terms of the Reidemeister torsion forms.

In this paper we analogously construct $\tau$ for $K$ in a general rational homology $3$-sphere and give $\mathcal{R}$ an orientation and a Riemannian metric induced by $\tau$.
Then we consider an involution $\iota$ on $\mathcal{R}$ which comes from the fact that the center of $SU_2$ is $\Z/2$ and show that $\iota$ is an orientation reversing isometry (Proposition \ref{prop_I}).
As a corollary we can detect the distribution of conjugacy classes of metabelian representations in $\mathcal{R}$ (Corollary \ref{cor_M}).
Next we assume that $E$ is irreducible and that $\partial E$ is incompressible and study the action on $\mathcal{R}$ by the outer automorphism group $\out(\pi_1 E)$.
More precisely, we consider the peripheral structure preserving outer automorphism group $\out_p(\pi_1 E)$ and show that the induced map $[\varphi]^*$ by $[\varphi] \in \out_p(\pi_1 E)$ is an isometry on an appropriate open subspace $\mathcal{R}_{\varphi} \subset \mathcal{R}$ and that the isometry preserves or reverses the orientation according to a sign $\delta_D^{[\varphi]}$ determined by each component $D$ of $\mathcal{R}_{\varphi}$ (Theorem \ref{thm_O}).
Here the basic ingredient is Waldhausen's theorem which allows us to take a homeomorphism of $E$ associated to $[\varphi]$. 
Finally, as a consequence of these results we introduce $2$-variable Reidemeister torsion $T_{D}(s, t)$ where $s$ is a coordinate on a component $D$ of $\mathcal{R}$ or $\mathcal{R}_{\varphi}$ and show the following symmetry of $T_{D}(s, t)$.
We take a meridional element $\mu \in \pi_1 E$ which is compatible with the orientation and the natural surjection $\alpha \colon \pi_1 E \to H_1(E) / \tor H_1(E) = \langle t \rangle$ which maps $\mu$ to $t$.
In the theorem the relation $\sim$ means that two functions are equal up to translation on $s$.

\begin{thm}[Theorem \ref{thm_S}]
(i) For any component $D$ of $\mathcal{R}$,
\[ T_{\iota(D)}(s, t) \sim -T_D(-s, -t). \]
(ii) Suppose that $E$ is irreducible and that $\partial E$ is incompressible.
Let $[\varphi] \in \out_p(\pi_1 E)$ such that $\alpha \circ \varphi(\mu) = t^{\pm 1}$ and $D$ a component of $\mathcal{R}_{[\varphi]}$.
Then
\[ T_{[\varphi]^*(D)}(s, t) \sim T_D(\delta_D^{[\varphi]} s, t^{\pm 1}). \]
\end{thm}

This paper is organized as follows.
In the next section we explain fundamental facts on $SU_2$-representation spaces of $\pi_1 E$ and on Reidemeister torsion of CW-complexes.
Section $3$ is devoted to construct the volume form $\tau$.
In Section $4$ we study the involution $\iota$ and the action of $\out_p(\pi_1 E)$.
In the last section we describe symmetry of the torsion function and have an example for the figure eight knot. \\

\section{Preliminaries}
In this paper all homology groups and cohomology groups are with respect to integral coefficients unless specifically noted.
Let $X$ be a connected CW-complex and $\rho \colon \pi_1 X \to GL_n(R)$ a linear representation over a commutative ring $R$.
We regard $R^n$ as a left $\Z[\pi_1 X]$-module by
\[ \gamma \cdot v :=  \rho(\gamma) v, \]
where $\gamma \in \pi_1 X$ and $v \in R^n$.
Then we define the twisted homology group and the twisted cohomology group associated to $\rho$ as follows:
\begin{align*}
H_i^{\rho}(X; R^n) &:= H_i(C_*(\widetilde{X}) \otimes_{\Z[\pi_1 X]} R^n), \\
H_{\rho}^i(X; R^n) &:= H^i(\hom_{\Z[\pi_1 X]}(C_*(\widetilde{X}), R^n)),
\end{align*}
where $\widetilde{X}$ is the universal covering of $X$.

\subsection{$SU_2$-representation spaces} \label{subsec_2_1}
We begin by briefly reviewing the $SU_2$-representation spaces of knot groups.

We fix a base point in $\partial E$ and take a longitude-meridian pair $\lambda$, $\mu \in \pi_1 \partial E$ which is compatible with the orientations.
By abuse of notation we use the same letters $\lambda$, $\mu$ for the images by the natural map $\pi_1 \partial E \to \pi_1 E$.
We regard $\hom(\pi_1 E, SU_2)$ as a topological space via the compact open topology, where $\pi_1 E$ carries the discrete topology and $SU_2$ the usual topology.
Let us denote by $\mathcal{R}$ the space of conjugacy classes of irreducible $SU_2$-representations such that 
\[ \dim H_{\Ad \circ \rho}^1(E; \mathfrak{su}_2) = 1 \quad \text{and} \quad \rho(\mu) \neq \pm I \]
for any $[\rho] \in \mathcal{R}$, where $\Ad \colon SU_2 \to \aut(\mathfrak{su}_2)$ is the adjoint representation and $I$ the identity matrix.
Note that a representation $\rho \colon \pi_1 E \to SU_2$ is irreducible if and only if $\rho$ is non-abelian.
It is easy to check that if the ambient space of $K$ is $S^3$, then no irreducible $SU_2$-representation satisfies that $\rho(\mu) = \pm I$.  

The representation spaces $\hom(\pi_1 E, SU_2)$ and $\hom(\pi_1 E, SL_2(\C))$ have usual real and complex affine algebraic structures.
We consider the involution $\sigma$ on $SL_2(\C)$ defined by
\[ \sigma(A) = (\overline{A}^T)^{-1}. \]
We also denote by $\sigma$ the induced involution on $\hom(\pi_1 E, SL_2(\C))$.
It can be easily seen that the fixed point set $\hom(\pi_1 E, SL_2(\C))^{\sigma}$ is nothing but $\hom(\pi_1 E, SU_2)$.
Following Weil \cite{We}, the Zariski tangent space $T_{\rho}Hom(\pi_1 E, SL_2(\C))$ can be identified with a subspace of the vector space $Z_{\Ad \circ \rho}^1(G; \mathfrak{sl}_2(\C))$ of group $1$-cocycles by the inclusion given by
\[ \left. \frac{d \rho_{t}}{dt} \right|_{t=0} ~\mapsto~ \left( \gamma \mapsto \left. \frac{d \rho_t (\gamma) \rho (\gamma^{-1})}{dt} \right|_{t=0} \right), \]
where $\rho_0 = \rho$ and $\gamma \in \pi_1 E.$
Thus we have the inclusion
\begin{equation} \label{eq_W}
T_{\rho}Hom(\pi_1 E, SU_2) = (T_{\rho}Hom(\pi_1 E, SL_2(\C)))^{\sigma_*} \hookrightarrow (Z_{\Ad \circ \rho}^1(G; \mathfrak{sl}_2(\C)))^{\sigma_*} = Z_{\Ad \circ \rho}^1(G; \mathfrak{su}_2).
\end{equation}
We refer the reader to \cite{P}.

The following proposition can be proved in the same way as in \cite[Proposition.~1]{HK}.
Therefore we omit the proof.
\begin{prop} \label{prop_R}
(i)The space $\mathcal{R}$ is a $1$-dimensional manifold. \\
(ii)The inclusion in \eqref{eq_W} induces an isomorphism
\[ T_{[\rho]} \mathcal{R} \cong H_{\Ad \circ \rho}^1(E; \mathfrak{su}_2). \]
\end{prop}

For a generator system $S = \{ \gamma_1, \dots. \gamma_m \}$ of $\pi_1 E$, we define a map $j_S \colon \mathcal{R} \to \R^{2^m - 1}$ by
\[ j_S([\rho]) = (\tr \rho(\gamma_{i_1} \dots \gamma_{i_k})), \]
where $1 \leq i_1 < \dots < i_k \leq m$.
We can also see that $j_S$ is an embedding for any $S$.
This follows from the fact that we can regard $\mathcal{R}$ as a subspace of the $SL_2(\C)$-character variety of $\pi_1 E$ and any point of $\mathcal{R}$ is smooth.
See \cite{CS} and \cite{GM} for the definition of the $SL_2(\C)$-character varieties of finitely generated groups.
In particular, for any $\gamma \in \pi_1 E$, the map $I_{\gamma} \colon \mathcal{R} \to \R$ defined by
\[ I_{\gamma}([\rho]) = \tr \rho(\gamma) \]
for $[\rho] \in \mathcal{R}$ is a smooth map and for any $[\rho] \in \mathcal{R}$, there exists $\gamma \in \pi_1 E$ such that $I_{\gamma}$ is an embedding in a neighborhood of $[\rho]$.

\subsection{Reidemeister torsion}
Next we review the definition of Reidemeister torsion, following Turaev \cite{T1}, \cite{T2}, and introduce normalized Reidemeister torsion for a knot exterior $E$.

For given bases $v$ and $w$ of a vector space, we define $[v / w]$ to be the determinant of the base change matrix from $w$ to $v$.

Let $C_* = (C_n \xrightarrow{\partial_n} C_{n-1} \to \cdots \to C_0)$ be a chain complex of finite dimensional vector spaces over a field $\F$.
For given bases $b_i$ of $\im \partial_{i+1}$ and $h_i$ of $H_i(C_*)$, we can choose a basis $b_i h_i b_{i-1}$ of $C_i$ as follows.
Choosing a lift of $h_i$ in $\ker \partial_i$ and combining it with $b_i$, we obtain a basis $b_i h_i$ of $\ker \partial_i$.
Then choosing a lift of $b_{i-1}$ in $C_i$ and combining it with $b_i h_i$, we obtain a basis $b_i h_i b_{i-1}$ of $C_i$.
\begin{defn}
For given bases $\boldsymbol{c} = \{ c_i \}$ of $C_*$ and $\boldsymbol{h} = \{ h_i \}$ of $H_*(C_*)$, we choose a basis $\{ b_i \}$ of $\im \partial_*$ and define
\[ \tau (C_*, \boldsymbol{c},\boldsymbol{h}) := (-1)^{|C_*|} \prod_{i=0}^n [b_i h_i b_{i-1} / c_i]^{(-1)^{i+1}} ~ \in \F^*, \]
where
\[ |C_*| := \sum_{j=0}^n (\sum_{i=0}^j \dim C_i)(\sum_{i=0}^j \dim H_i(C_*)). \]
\end{defn}
It can be easily checked that $\tau (C_*,\boldsymbol{c},\boldsymbol{h})$ does not depend on the choices of $b_i$ and $b_i h_i b_{i-1}$.

Let $X$ be a connected finite CW-complex with cells $\{ e_i \}$.
For two lifts $\{ \tilde{e}_i \}$ and $\{ \tilde{e}_i' \}$ in the universal covering $\widetilde{X}$, we set
\[ \{ \tilde{e}_i' \} / \{ \tilde{e}_i \} := \sum_i (-1)^{\dim e_i} \tilde{e}_i' / \tilde{e}_i ~ \in H_1(X), \]
where $\tilde{e}_i' / \tilde{e}_i$ is the element $h$ of $H_1(X)$ such that $\tilde{e}_i' = \tilde{e}_i \cdot h$.
Two lifts $\{ \tilde{e}_i \}$ and $\{ \tilde{e}_i' \}$ are called equivalent if $\{ \tilde{e}_i' \} / \{ \tilde{e}_i \} = 0$.
An equivalence class is called an \textit{Euler structure} of $X$ and the set of Euler structures is denoted by $\eul(X)$.
For $h \in H_1(X)$ and $[\{ \tilde{e}_i \}] \in \eul(X)$, let $[\{ \tilde{e}_i \}] \cdot h$ be the class of any lift $\{ \tilde{e}_i \}$ with $\{ \tilde{e}_i' \} / \{ \tilde{e}_i \} = h$.
This defines a free and transitive action of $H_1(X)$ on $\eul(X)$.
It easily follows that for any subdivision $X'$ of $X$, there is a canonical $H_1(X)$-bijection between $\eul(X)$ and $\eul(X')$.
By a \textit{homology orientation} of $X$ we mean an orientation of the homology group $H_*(X; \R) = \bigoplus_i H_i(X; \R)$ as a vector space.
\begin{defn}
For a representation $\rho \colon \pi_1 X \to GL_n(\F)$, an Euler structure $\mathfrak{e}$ and a homology orientation $\mathfrak{o}$, we define the \textit{Reidemeister torsion} $T_{\rho}(X, \mathfrak{e}, \mathfrak{o})$ associated to $\rho$, $\mathfrak{e}$ and $\mathfrak{o}$ as follows.
We choose a lift $\{ \tilde{e}_i \}$ of cells $\{ e_i \}$ representing $\mathfrak{e}$ and bases $\boldsymbol{h}$ of $H_*(X; \R)$ which is positively oriented with respect to $\mathfrak{o}$ and $\langle f_1, \dots, f_n \rangle$ of $\F^n$.
Then we define
\[ T_{\rho}(X, \mathfrak{e}, \mathfrak{o}) :=
\begin{cases}
\tau_0^n \tau (C_*^{\rho}(X; \F^n), \tilde{\boldsymbol{c}}) &\text{if $H_*^{\rho}(X; F^n)$ vanishes}, \\
0 &\text{otherwise,}
\end{cases} \]
where
\begin{align*}
\tau_0 &:= \sgn \tau (C_*(X;\R), \boldsymbol{c}, \boldsymbol{h}), \\
\boldsymbol{c} &:= \langle e_1, \dots , e_{dim C_*(X)} \rangle, \\
\tilde{\boldsymbol{c}} &:= \langle \tilde{e}_1 \otimes f_1, \dots , \tilde{e}_1 \otimes f_n, \dots , \tilde{e}_{dim C_*(X)} \otimes f_1, \dots , \tilde{e}_{dim C_*(X)} \otimes f_n \rangle.
\end{align*}
\end{defn}
It is known that $T_{\rho}(X, \mathfrak{e}, \mathfrak{o})$ does not depend on the choices of $\tilde{e}_i$, $\boldsymbol{h}$ and $\langle f_1, \dots, f_n \rangle$ and $T_{\rho}(X, \mathfrak{e}, \mathfrak{o})$ is invariant under cellular subdivision.
\begin{rem}
For a link exterior of $S^3$, given a presentation of the link group, Reidemeister torsion can be computed efficiently using Fox calculus (cf.\ e.\ g.\ \cite{KL}, \cite{Ki}).
\end{rem}

Turaev defined an involution $\mathfrak{e} \mapsto \mathfrak{e}^{-1}$ on $\eul(X)$ and a map $c \colon \eul(X) \to H_1(X)$ which maps $\mathfrak{e}$ to $\mathfrak{e} / \mathfrak{e}^{-1}$ for a compact $3$-manifold $X$ whose boundary consists of tori (also for a closed odd-dimensional manifold $X$).
The map $c$ satisfies that
\[ c(\mathfrak{e} \cdot h) = c(\mathfrak{e}) + 2 h \]
for any $h \in H_1(X)$.
We identify an infinite cyclic group $\langle t \rangle$ with $H_1(E) / \tor H_1(E)$.
Let $\alpha \colon \pi_1 E \to \langle t \rangle$ be the map which maps $\mu$ to $t$ and $k \colon \eul(E) \to \Z$ the map defined by
\begin{equation} \label{eq_c}
[c(\mathfrak{e})] = t^{k(\mathfrak{e})}.
\end{equation}
It is known that $k(\mathfrak{e})$ is an odd number.
(See \cite{T2}.)

\begin{defn}
Let $\rho \colon \pi_1 E \to SL_n(\F)$ be a representation.
We choose $\mathfrak{e} \in \eul(E)$ and then define
\[ \widetilde{T}_{\rho}(t) := t^{-\frac{n k(\mathfrak{e})}{2}} T_{\alpha \otimes \rho}(E, \mathfrak{e}, \{ [pt], [\mu] \}) ~\in \F(t^{\frac{1}{2}}), \]
where $\alpha \otimes \rho \colon \pi_1 E \to GL_n(\F(t))$ is a representation which maps $\gamma \in \pi_1 E$ to $\alpha(\gamma) \rho(\gamma)$.
We call it the \textit{normalized Reidemeister torsion} associated to $\rho$.
\end{defn}

\begin{lem}
The normalized Reidemeister torsion $\widetilde{T}_{\rho}$ does not depend on the choice of $\mathfrak{e}$.
\end{lem}

\begin{proof}
From \eqref{eq_c} and the definitions, for any $h \in H_1(E)$ with $[h] = t^m$, we have
\begin{align*}
k(\mathfrak{e} \cdot h) &= k(\mathfrak{e}) + 2m, \\
T_{\alpha \otimes \rho}(E, \mathfrak{e} \cdot h, \{ [pt], [\mu] \}) &= t^{m n} T_{\alpha \otimes \rho}(E, \mathfrak{e}, \{ [pt], [\mu] \}),
\end{align*}
which proves the lemma.
\end{proof}
\section{Torsion volume forms}
In this section we construct a volume form $\tau$ on $\mathcal{R}$ via non-acyclic Reidemeister torsion, slightly generalizing Dubois' form in \cite{D1}.

\begin{lem} \label{lem_v}
For any $[\rho] \in \mathcal{R}$,
\begin{align*}
\dim H_{\Ad \circ \rho}^0(E; \mathfrak{su}_2) &= 0, \\
\dim H_{\Ad \circ \rho}^2(E; \mathfrak{su}_2) &= 1.
\end{align*}
\end{lem}

\begin{proof}
Since $\rho$ is non-abelian,
\[ H_{\Ad \circ \rho}^0(E; \mathfrak{su}_2) = \mathfrak{su}_2^{\Ad \circ \rho(\pi_1 E)} = 0. \]
Considering
\[ \sum_{i = 0}^2 (-1)^i \dim H_{\Ad \circ \rho}^i(E; \mathfrak{su}_2) = 3 \chi(E) = 0, \]
we have
\[ \dim H_{\Ad \circ \rho}^2(E; \mathfrak{su}_2) = \dim H_{\Ad \circ \rho}^1(E; \mathfrak{su}_2) = 1. \]
\end{proof}

The Killing form of $\mathfrak{su}_2$ induces non-degenerate cup products
\begin{align}
\cup \colon &H_{\Ad \circ \rho}^q(E; \mathfrak{su}_2) \times H_{\Ad \circ \rho}^{3-q}(E, \partial E; \mathfrak{su}_2) \to H_{\Ad \circ \rho}^3(E, \partial E; \mathfrak{su}_2), \label{eq_c1} \\
&H_{\Ad \circ \rho}^q(\partial E; \mathfrak{su}_2) \times H_{\Ad \circ \rho}^{2-q}(\partial E; \mathfrak{su}_2) \to H_{\Ad \circ \rho}^2(\partial E; \mathfrak{su}_2). \label{eq_c2}
\end{align}

Let $\psi$ be the map
\[ H_{\Ad \circ \rho}^2(E; \mathfrak{su}_2) \to H_{\Ad \circ \rho}^2(\partial E; \mathfrak{su}_2) \to H_{\Ad \circ \rho}^0(\partial E; \mathfrak{su}_2)^*, \]
where these maps are the induced homomorphisms by the natural inclusion and \eqref{eq_c2} respectively.
For $[\rho] \in \mathcal{R}$, we can take unique pairs $(\theta_{\rho}, P_{\rho}) \in (0, \pi) \times \mathfrak{su}_2$ such that
\[ \rho(\mu) = I \cos \theta + P_{\rho} \sin \theta_{\rho}, \]
where we regard $\mathfrak{su}_2$ as the set of skew-symmetric trace free matrices.
Notice that
\[ P_{A \rho A^{-1}} = A P_{\rho} A^{-1} \]
for $A \in SU_2$.

\begin{lem}
For any $[\rho] \in \mathcal{R}$, $\psi$ is an isomorphism and
\[ H_{\Ad \circ \rho}^0(\partial E; \mathfrak{su}_2) = \R P_{\rho}. \]
\end{lem}

\begin{proof}
It is easily seen that
\[ (\Ad \circ \rho(\mu))(P_{\rho}) = (\Ad \circ \rho(\lambda))(P_{\rho}) = P_{\rho}, \]
and so we have
\[ H_{\Ad \circ \rho}^0(\partial E; \mathfrak{su}_2) = \mathfrak{su}_2^{\Ad \circ \rho(\pi_1 \partial E)} = \R P_{\rho}, \]
where
\[ \mathfrak{su}_2^{\Ad \circ \rho(\pi_1 E)} = \{ \xi \in \mathfrak{su}_2 ~;~ \Ad \circ \rho(\gamma)(\xi) = \xi \text{ for any } \gamma \in \pi_1 E \}. \]

It remains to prove that the homomorphism $H_{\Ad \circ \rho}^2(E; \mathfrak{su}_2) \to H_{\Ad \circ \rho}^2(\partial E; \mathfrak{su}_2)$ is an isomorphism.
From Lemma \ref{lem_v} and \eqref{eq_c1}
\[ \dim H_{\Ad \circ \rho}^3(E, \partial E; \mathfrak{su}_2) = \dim H_{\Ad \circ \rho}^0(E; \mathfrak{su}_2) = 0. \]
Therefore from the long exact sequence for the pair ($E$, $\partial E$) the above homomorphism is surjective.
From \eqref{eq_c1}
\[ \dim H_{\Ad \circ \rho}^2(\partial E; \mathfrak{su}_2) = \dim H_{\Ad \circ \rho}^0(\partial E; \mathfrak{su}_2) = 1 \]
and from Lemma \ref{lem_v} $\dim H_{\Ad \circ \rho}^2(E; \mathfrak{su}_2)$ is also $1$, which deduces the desired conclusion.
\end{proof}

For $[\rho] \in \mathcal{R}$, we set 
\[ h_{\rho} = \psi^{-1}(P_{\rho}). \]
\begin{defn} \label{def_V}
We choose a lift $\tilde{e}_i$ in $\widetilde{E}_K$ of each cell $e_i$ in $E_K$ and a basis $\langle \xi_1, \xi_2, \xi_3 \rangle$ of $\mathfrak{su}_2$.
At $[\rho] \in \mathcal{R}$ a linear form $\tau_{[\rho]} : T_{[\rho]} \mathcal{R} \to \R$ is defined by
\[ \tau_{[\rho]}(v) :=
\begin{cases}
\tau_0 \tau(C_{\Ad \circ \rho}^{-*}(E; \mathfrak{su}_2), \tilde{\boldsymbol{c}}, \langle v, h_{\rho} \rangle ) &\text{if } v \neq 0, \\
0 &\text{if } v = 0, 
\end{cases}
\]
where
\begin{align*}
\tau_0 &:= \sgn \tau (C^{-*}(E; \R),\boldsymbol{c} ,\langle [pt]^*, [\mu]^* \rangle), \\
\boldsymbol{c} &:= \langle e_1^*, \dots , e_{dim C_*}^* \rangle, \\
\tilde{\boldsymbol{c}} &:= \langle \tilde{e}_{1, 1}, \tilde{e}_{1, 2}, \tilde{e}_{1, 3}, \dots , \tilde{e}_{dim C_*(E), 1}, \tilde{e}_{dim C_*(E), 2}, \tilde{e}_{dim C_*(E), 3} \rangle
\end{align*}
and $\tilde{e}_{i, j}$ is a cochain which maps $\tilde{e}_i$ to $\xi_j$ and $\tilde{e}_{i'}$ to $0$ for $i' \neq i$.
\end{defn}

It can be checked that the linear form $\tau$ does not depend on the choices of $\{ \tilde{e}_i \}$, $\langle \xi_1, \xi_2, \xi_3 \rangle$ and $\rho$, and so $\tau$ is well defined as a volume form on $\mathcal{R}$ as Dubois' form is.

\section{Isometries}
\subsection{$(-1)$-involution}
We give $\mathcal{R}$ the orientation and the Riemannian metric induced by $\tau$.

There is an involution $\iota$ on $\hom(\pi_1 E, SU_2)$ defined by
\[ \iota(\rho)(\gamma) = \alpha(\gamma) |_{t=-1} \rho(\gamma), \]
for $\gamma \in \pi_1 E$.
Recall that $\alpha$ is the map $\pi_1 E \to H_1(E) / \tor H_1(E) = \langle t \rangle.$
We write $\rho^* := \iota(\rho)$ for simplicity.
\begin{lem} \label{lem_i}
The involution $\iota$ induces a diffeomorphism $\mathcal{R} \to \mathcal{R}$.
\end{lem}

\begin{proof}
Since we have $\Ad \circ \rho^* = \Ad \circ \rho$ and $\rho^*(\mu) = -\rho(\mu)$, $\iota$ naturally induces an involution on $\mathcal{R}$.

From the observation in Section \ref{subsec_2_1}, for any point $[\rho] \in \mathcal{R}$, there exist $\gamma$, $\gamma' \in \pi_1 E$ and open neighborhoods $U$, $U'$ of $[\rho]$, $[\rho^*]$ respectively such that $\iota(U) \subset U'$ and $I_{\gamma} |_U$, $I_{\gamma'} |_{U'}$ are embeddings.
Since 
\[ I_{\gamma'} |_{U'} \circ \iota \circ (I_{\gamma} |_U)^{-1} = \alpha(\gamma') |_{-1} I_{\gamma'} \circ (I_{\gamma} |_U)^{-1}, \]
the induced map is smooth, which proves the lemma.
\end{proof}

By abuse of notation we also denote by $\iota$ the above induced map.
\begin{prop} \label{prop_I}
The involution $\iota$ on $\mathcal{R}$ is an orientation reversing isometry, i.\ e.\ ,
\[ \iota^* \tau = - \tau. \]
\end{prop}

\begin{proof}

For any $v \in T_{[\rho]} \mathcal{R}$, choose a smooth family $\{ [\rho_t] \} \in \mathcal{R}$ such that $v = \left. \frac{d [\rho_t]}{dt} \right|_{t=0}$, where $\rho_0 = \rho$.
Since
\[ \left. \frac{d}{dt} \rho_t^*(\gamma) \rho^*(\gamma^{-1}) \right|_{t=0} = \left. \frac{d}{dt} \rho_t(\gamma) \rho(\gamma^{-1}) \right|_{t=0}, \]
we obtain $\iota_*(v) = v$ by Proposition \ref{prop_R} (ii).
Considering $P_{\rho^*} = -P_{\rho}$, we have $h_{\rho^*} = -h_{\rho}$.
These give
\[ \tau_{[\rho^*]}(\iota_*(v)) = -\tau_{[\rho]}(v). \]
\end{proof}

A representation $\rho$ of a group $G$ is called metabelian if $\rho([G, G])$ is abelian, where $[G, G]$ is the commutator subgroup of $G$.
We denote by $\mathcal{S}$ the set of conjugacy classes of irreducible metabelian $SU_2$-representations of $\pi_1 E$.
\begin{prop}[\cite{NY}]
Let $K$ be a knot in an integral homology $3$-sphere.
Then the fixed point set of $\iota$ in the set of conjugacy classes of irreducible $SU_2$-representations is equal to $\mathcal{S}$.
\end{prop}

Since an orientation reversing isometry on an arc or a circle which fixes a point is uniquely determined, we immediately deduce the following corollary.
\begin{cor} \label{cor_M}
(i)If $[\rho] \in \mathcal{S}$ is a point on an arc component $A$ of $\mathcal{R}$, $[\rho]$ is a unique point on $\mathcal{S} \cap A$ and $[\rho]$ is at the center of $A$. \\
(ii)If $[\rho] \in \mathcal{S}$ is a point on a circle component $C$ of $\mathcal{R}$, there are exactly two points on $\mathcal{S} \cap C$ and they are antipodal to each other.
\end{cor}

\begin{rem}
In \cite[Proposition 5.\ 3.]{BLZ} Boyer, Luft and Zhang show that if $E$ is hyperbolic and fibered over $S^1$, then $\mathcal{S} \subset \mathcal{R}$.
Furthermore, it is implicit from \cite{HK} that if the ambient space of $K$ is $S^3$ and the knot determinant $|\Delta_K(-1)|$ is a prime number, then $\mathcal{S} \subset \mathcal{R}$.
\end{rem}

\subsection{Automorphism group actions}
In this subsection we assume that $E$ is irreducible and that $\partial E$ is incompressible.

The automorphism group $\aut(\pi_1 E)$ acts on $\hom(\pi_1 E, SU_2)$ via pullback.
We write $\varphi^* \rho := \rho \circ \varphi$ for $\varphi \in \aut(\pi_1 E)$ and $\rho \in \hom(\pi_1 E, SU_2)$.
This action naturally induces the action of $\out(\pi_1 E)$ on the set of conjugacy classes of representations.
For $[\varphi] \in \out(\pi_1 E)$, we denote by $\mathcal{R}_{[\varphi]}$ the set of $[\rho] \in \mathcal{R}$ which satisfies that $\varphi^* \rho(\mu) \neq \pm I$.
Note that if $K$ is a knot in $S^3$, then $\mathcal{R}_{[\varphi]} = \mathcal{R}$ for any $[\varphi] \in \out(\pi_1 E)$.
\begin{lem}
An outer automorphism $[\varphi]$ induces a diffeomorphism $[\varphi]^* \colon \mathcal{R}_{[\varphi]} \to \mathcal{R}_{[\varphi^{-1}]}$.
\end{lem}

\begin{proof}
We can proceed analogously to the proof of Lemma \ref{lem_i}.

Since there is an isomorphism $\varphi^* \colon H_{\Ad \circ \rho}^*(\pi_1 E; \mathfrak{su}_2) \to H_{\Ad \circ \varphi^* \rho}^*(\pi_1 E; \mathfrak{su}_2)$ induced by $\varphi$, $[\varphi]$ induces a map $\mathcal{R}_{[\varphi]} \to \mathcal{R}_{[\varphi^{-1}]}$.

For any point $[\rho] \in \mathcal{R}_{[\varphi]}$, there exist $\gamma$, $\gamma' \in \pi_1 E$ and open neighborhoods $U \subset \mathcal{R}_{[\varphi]}$, $U' \subset \mathcal{R}_{[\varphi^{-1}]}$ of $[\rho]$, $[\varphi^* \rho]$ respectively such that $[\varphi]^*(U) \subset U'$ and $I_{\gamma} |_U$, $I_{\gamma'} |_{U'}$ are embeddings.
Since 
\[ I_{\gamma'} |_{U'} \circ [\varphi]^* \circ (I_{\gamma} |_U)^{-1} = I_{\varphi(\gamma')} \circ (I_{\gamma} |_U)^{-1}, \]
the induced map is smooth, and the lemma follows.
\end{proof}

We denote by $\out_p(\pi_1 E)$ the set of outer automorphisms which preserve the peripheral structure, i.\ e.,
\[ \varphi(\langle \lambda, \mu \rangle) = \langle \lambda, \mu \rangle \]
for $[\varphi] \in \out_p(\pi_1 E)$.
For $[\varphi] \in \out_p(\pi_1 E)$, since $\lambda$ is a generator of $\ker(\langle \lambda, \mu \rangle)$, there exists $\delta \in \{ \pm 1 \}$ such that
\[ \varphi(\lambda) = \lambda^{\delta}.\]
Moreover, for $[\rho] \in \mathcal{R}_{[\varphi]}$, since $[\rho(\mu), \varphi^* \rho(\mu)] = 1$, there exist $\delta' \in \{ \pm 1 \}$ and $\theta' \in (0, \pi)$ such that
\[ \varphi^* \rho(\mu) = I \cos \theta' + \delta' P_{\rho} \sin \theta',\]
and so
\[ P_{\varphi^* \rho} = \delta' P_{\rho}. \]
\begin{defn} \label{def_D}
For $[\varphi] \in \out_p(\pi_1 E)$ and a component $D \in \mathcal{R}_{[\varphi]}$, choosing $[\rho] \in D$, we define
\[ \delta_D^{[\varphi]} := \delta \delta'. \]
\end{defn}

An easy verification shows that the sign $\delta_D^{[\varphi]}$ depends only on $[\varphi]$ and $D$.

\begin{thm} \label{thm_O}
For $[\varphi] \in \out_p(\pi_1 E)$ and a component $D \in \mathcal{R}_{[\varphi]}$, the following equality holds:
\[ [\varphi]^* \tau |_D = \delta_D^{[\varphi]} \tau |_{[\varphi]^*(D)}. \]
\end{thm}

\begin{proof}
By Waldhausen's theorem there is a homeomorphism $f$ which preserves the base point and satisfies that $f_* = \varphi$.
We take lifts $\tilde{f} \colon \widetilde{E} \to \widetilde{E}$ and $\tilde{e}_i \subset \widetilde{E}$ for each cell $e_i \subset E$.
Then $f$ induces an isomorphism $f^* \colon C_{\Ad \circ \rho}^*(f(E); \mathfrak{su}_2) \to C_{\Ad \circ \varphi^* \rho}^*(E; \mathfrak{su}_2)$ defined by
\[ (\tilde{f}(\tilde{e}_i) \mapsto \xi) \mapsto (\tilde{e} \mapsto \xi). \]
We have the following commutative diagram:
\[
\begin{CD}
T_{[\rho]}\mathcal{R} @= H_{\Ad \circ \rho}^1(\pi_1 E;\mathfrak{su}_2) @= H_{\Ad \circ \rho}^1(E;\mathfrak{su}_2) \\
@V([\varphi]^*)_*VV @V\varphi^*VV @Vf^*VV \\
T_{[\varphi^* \rho]}\mathcal{R} @= H_{\Ad \circ \varphi^* \rho}^1(\pi_1 E;\mathfrak{su}_2) @= H_{\Ad \circ \varphi^* \rho}^1(E;\mathfrak{su}_2)
\end{CD}
\]

Let $[\rho]$ be a point in $D$.
We use the notation of Definition \ref{def_V} and set
\begin{align*}
\boldsymbol{c}_f &:= \langle f(e_1)^*, \dots , f(e_{dim C_*})^* \rangle, \\
\tilde{\boldsymbol{c}}_f &:= \langle \tilde{f}(\tilde{e}_1)_1, \tilde{f}(\tilde{e}_1)_2, \tilde{f}(\tilde{e}_1)_3, \dots , \tilde{f}(\tilde{e}_{dim C_*(E)})_1, \tilde{f}(\tilde{e}_{dim C_*(E)})_2, \tilde{f}(\tilde{e}_{dim C_*(E)})_3 \rangle,
\end{align*}
where $\tilde{f}(\tilde{e}_i)_j$ is a cochain in $C_{\Ad \circ \rho}^*(f(E);\mathfrak{su}_2)$ which maps $\tilde{f}(\tilde{e}_i)$ to $\xi_j$ and $\tilde{f}(\tilde{e}_{i'})$ to $0$ for $i' \neq i$.
Using isomorphisms induced by $f$, we have
\begin{align*}
\tau_{[\rho]}(v) &= (\sgn \tau (C^{-*}(f(E); \R),\boldsymbol{c}_f ,\langle [pt]^*, [\mu]^* \rangle)) \tau(C_{\Ad \circ \rho}^{-*}(E; \mathfrak{su}_2), \tilde{\boldsymbol{c}}_f, \langle v, h_{\rho} \rangle ) \\
&= [[\varphi(\mu)]/ [\mu]] \tau_0 \cdot [h_{\varphi^* \rho}/ f^*(h_{\rho})] \tau(C_{\Ad \circ \varphi^* \rho}^{-*}(E; \mathfrak{su}_2), \tilde{\boldsymbol{c}}, \langle ([\varphi]^*)_*(v), h_{\varphi^* \rho} \rangle ) \\
&= \delta_D^{[\varphi]} \tau_{[\varphi^* \rho]}(([\varphi]^*)_*(v)).
\end{align*}
\end{proof}

In the case where the ambient space is $S^3$ the following lemma is useful. 
\begin{lem}
Suppose $K$ is a knot in $S^3$.
Then for $[\varphi] \in \out_p(\pi_1 E)$ and a component $D \in \mathcal{R}$,
\begin{align*}
\varphi(\mu) &= \mu^{\pm 1}, \\
\delta_D^{[\varphi]} &= \det \varphi |_{\langle \lambda, \mu \rangle.}
\end{align*}
\end{lem}

\begin{proof}
We can write
\[ \varphi(\mu) = \lambda^l \mu^{\delta''},\]
where $l \in \Z$ and $\delta'' = \pm 1$.
We can again consider the homeomorphism $f$ in the proof of Theorem \ref{thm_O}.
Therefore the result of $\frac{\delta''}{l}$-surgery along $K$ is homeomorphic to $S^3$.
But by the theorem of Gordon and Luecke \cite{GL} we have $l = 0$.
Now it is easily seen that
\[ P_{\varphi^* \rho} = \delta'' P_{\rho}, \]
which establishes the second equality.
\end{proof}

\section{Symmetry of torsion functions}
Finally we apply the study in the previous sections to investigate symmetry of the torsion function on $\mathcal{R}$.

\begin{defn} \label{def_T}
Let $D$ be a component of $\mathcal{R}$ or $\mathcal{R}_{[\varphi]}$ for $[\varphi] \in \out(\pi_1 E)$.
We fix $[\rho_0] \in D$ and choose $[\rho_s] \in D$ for an appropriate $s \in \R$ so that
\[ \int_0^1 p^* \tau dt = s \]
for some smooth path $p \colon [0, 1] \to D$ from $[\rho_0]$ to $[\rho_s]$.
Then we define
\[ T_D(s, t) := \widetilde{T}_{\rho_s}(t). \]
\end{defn}
It is clear that the function $T_D$ is well defined up to translation on $s$ and that for a circle component $D$, $T_D$ is a periodic function on $s$ whose period is $\vol D$.

For two functions $f(s)$ and $g(s)$ on open intervals $I$ and $J \subset \R$ respectively, we write $f(s) \sim g(s)$ if $I$ and $J$ have the same length and there exists $s_0 \in \R$ such that $f(s) = g(s + s_0)$ for any $s \in I$.

\begin{thm} \label{thm_S}
(i) For any component $D$ of $\mathcal{R}$,
\[ T_{\iota(D)}(s, t) \sim -T_D(-s, -t). \]
(ii) Suppose that $E$ is irreducible and that $\partial E$ is incompressible.
Let $[\varphi] \in \out_p(\pi_1 E)$ such that $\alpha \circ \varphi(\mu) = t^{\pm 1}$ and $D$ a component of $\mathcal{R}_{[\varphi]}$.
Then
\[ T_{[\varphi]^*(D)}(s, t) \sim T_D(\delta_D^{[\varphi]} s, t^{\pm 1}). \]
\end{thm}

This theorem is a simple corollary of Proposition \ref{prop_I}, Theorem \ref{thm_O} and the following lemma.
\begin{lem}
(i) For $\rho \in \hom(\pi_1 E, SU_2)$,
\[ \widetilde{T}_{\rho^*}(t) = -\widetilde{T}_{\rho}(-t). \]
(ii) Suppose that $E$ is irreducible and that $\partial E$ is incompressible.
Let $\varphi$ be a peripheral structure preserving automorphism such that $\alpha \circ \varphi(\mu) = t^{\pm 1}$.
Then for $\rho \in \hom(\pi_1 E, SU_2)$,
\[ \widetilde{T}_{\varphi^* \rho}(t) = \widetilde{T}_{\rho}(t^{\pm 1}). \]
\end{lem}

\begin{proof}
(i) is straightforward from the definition and the fact that $k(\mathfrak{e})$ is odd for any $\mathfrak{e} \in \eul(E)$.

Suppose the assumptions of (ii).
By Waldhausen's theorem we can take a base point preserving homeomorphism $f \colon E \to E$ such that $f_* = \varphi$ and then $f$ induces an isomorphism $C_*^{(\alpha \otimes \rho) \circ \varphi}(E; \C(t)^2) \to C_*^{\alpha \otimes \rho}(f(E); \C(t)^2)$.
Computing torsion of these chain complexes, we obtain the desired equality in (ii).
\end{proof}

\begin{exmp}
\begin{figure} \label{fig_F}
\centering
\includegraphics[width=5cm, clip]{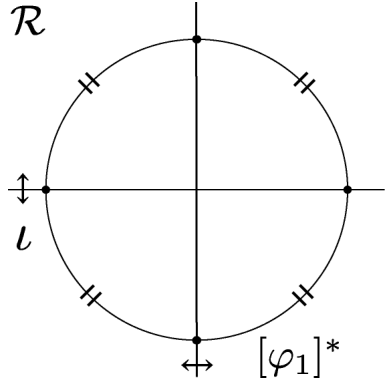}
\caption{The space $\mathcal{R}$ of $4_1$}
\end{figure}

Let $K$ be the figure eight knot $4_1$ in $S^3$.
A presentation of $\pi_1 E$ is given by
\[ \pi_1 E = \langle x, y \mid [x^{-1}, y] x = y [x^{-1}, y] \rangle, \]
where we can set
\begin{align*}
\lambda &= y x^{-1} y^{-1} x^2 y^{-1} x^{-1} y, \\
\mu &= x.
\end{align*}
From the result of Klassen in \cite{Kl} the space of conjugacy classes of irreducible $SU_2$-representations is a circle.
By direct computations we can check that $\mathcal{R}$ coincides with the whole space and
\[ \out(\pi_1 E) = \langle [\varphi_1], [\varphi_2] \mid [\varphi_1]^4 = [\varphi_2]^2 = ([\varphi_1] [\varphi_2])^2 = 1 \rangle, \]
where $\varphi_1$ is induced by an amphicheiral map and defined by
\begin{align*}
x &\mapsto x^{-1}, \\
y &\mapsto y x^{-1} y^{-1}
\end{align*}
and $\varphi_2$ is induced by an inversion map and defined by
\begin{align*}
x &\mapsto x^{-1}, \\
y &\mapsto y^{-1}.
\end{align*}
We can check that $[\varphi_1]^*$ is the orientation reversing isometry whose symmetric axis is orthogonal to that of $\iota$ as in Figure \ref{fig_F} and $[\varphi_2]^*$ is the identity map.
Computations give
\[ T_{\mathcal{R}}(s, t) = t - 2 f(s) + t^{-1}, \]
where $f(s)$ is a periodic function described in Figure \ref{fig_S}.
Using the notation in Definition \ref{def_T}, we obtain
\[ f(s) = \tr \rho_s(\mu). \]

\begin{figure} \label{fig_S}
\centering
\includegraphics[width=5cm, clip]{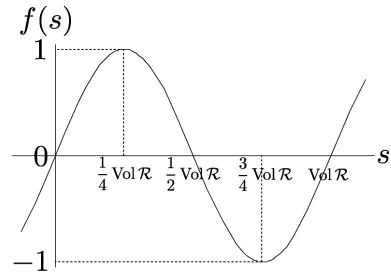}
\caption{The function $f(s)$}
\end{figure}
\end{exmp}

\noindent
\textbf{Acknowledgement.}
The author would like to express his gratitude to Toshitake Kohno for his encouragement and helpful suggestions.
He also would like to thank Hiroshi Goda, Teruaki Kitano, Takayuki Morifuji and Yoshikazu Yamaguchi for fruitful discussions and advices.
This research is supported by JSPS Research Fellowships for Young Scientists.


\end{document}